\documentclass{article}

\usepackage{amsmath,amssymb}
\usepackage[latin1]{inputenc}
\usepackage[dvips]{graphics}
\input{xy}
\xyoption{poly}
\xyoption{2cell}
\xyoption{all}

\newcommand{\qed}{$\ \hfill\Box $\bigskip}

\newcommand{\za}{\alpha}

\newcommand{\zg}{\gamma}

\newcommand{\zS}{\Sigma}

\newtheorem{thm}{Theorem}[section]

\newtheorem{cor}[thm]{Corollary}
\newtheorem{lem}[thm]{Lemma}
\newtheorem{example}[thm]{Example}
\newtheorem{definition}{Definition} 
\newtheorem{rem}[thm]{Remark}

\newenvironment{pf}{{Proof}.}

\begin{document}
\title{A cluster expansion formula ($A_n$ case)}
\author{R. Schiffler}
\date{}
\maketitle


\begin{abstract}
We consider the Ptolemy cluster algebras, which are  cluster algebras
of finite type $A$ (with non-trivial coefficients) that have been
described by Fomin and Zelevinsky using triangulations of a regular polygon.
Given any seed $\zS$ in a Ptolemy cluster algebra,
we present a formula for the expansion of an arbitrary cluster
variable in terms of the cluster variables of the
seed  $\zS$.
Our formula is given in a combinatorial way, using paths on a
triangulation of the polygon that
corresponds to the seed $\zS$.
 \end{abstract}
\maketitle

\setcounter{section}{0}
\addtocounter{section}{-1}


\begin{section}{Introduction}\label{section intro}
Ptolemy cluster algebras have been introduced by Fomin and
Zelevinsky in \cite[section 12]{FZ2} as  examples of  cluster algebras
of type $A$. 
The Ptolemy algebra of rank $n$ is described using 
the triangulations
of a regular polygon with $n+3 $ vertices. In this description, the
seeds of the cluster algebra are in bijection with the
triangulations of the polygon. The cluster of the seed corresponds to
the diagonals, while the coefficients of the seed correspond to the
boundary edges of the triangulation. 
The Laurent phenomenon \cite{FZ1} states that, given an arbitrary
seed $\zS$ one can write any cluster variable of 
the cluster algebra as a Laurent polynomial in the cluster variables
and the coefficients of the seed $\zS$. 

 The main result of this paper is an explicit formula for these
 Laurent polynomials, see Theorem
 \ref{thm 1}. Each term of
 the Laurent  polynomial is given by a path on the
 triangulation corresponding to the seed $\zS$.
As an application, we prove the positivity conjecture of Fomin and
Zelevinsky \cite{FZ1} for Ptolemy algebras, see
Corollary \ref{cor 1}. 

There is an interesting connection between our work and that of Propp
\cite{P} who uses perfect matchings arising from the triangulation 
to calculate these Laurent polynomials.

The polygon model has also been used in \cite{CCS1} to construct
the cluster category associated to the cluster algebra, compare \cite{BMRRT}. 
In that context, the cluster algebra has trivial coefficients and 
our formula naturally applies to that situation, see Remark \ref{rem 3}.
Therefore, there is
 an interesting intersection with the work of
Caldero and Chapoton \cite{CC},
who have obtained a formula for cluster expansions when the given seed
 is  acyclic, 
meaning that each triangle in the triangulation has at least one side
on the boundary of the polygon.
Their formula, and its generalization by Caldero and Keller
\cite{CK2}, also applies to cluster algebras (with trivial
 coefficients) of other types, but, again, only in the case where the
 seed $\zS$ is acyclic. Their description uses the representation theory
of finite dimensional algebras and is very different from ours.

We have been informed that our cluster expansion formula was also known to  
Carroll and Price  \cite{CP} and  to Fomin and Zelevinsky \cite{FZ3}.

The author thanks Hugh Thomas for interesting discussions on the subject.

\end{section} 



\begin{section}{Cluster expansions in the Ptolemy algebra}\label{sect 1}

Throughout this paper,
let $n$ be a positive integer and let $P$ be a regular  polygon with
$n+3$ vertices. 
A \emph{diagonal} in $P$ is a line segment connecting two non-adjacent
vertices of $P$ and two diagonals are said to be \emph{crossing} if they
intersect in the interior of the polygon.
A \emph{triangulation} $T$ is a maximal set of non-crossing
diagonals together with all boundary edges. Any triangulation $T$ has
$n$ diagonals and $n+3$ boundary edges.
Denote the boundary edges of $P$ by $T_{n+1},\ldots,T_{2n+3}$.

\begin{subsection}{The Ptolemy cluster algebra }\label{sect ptolemy}
We recall some facts about the Ptolemy cluster algebra of rank $n$ from
\cite[section 12.2]{FZ2}. The cluster variables $x_M$ of this algebra are in
bijection with the diagonals $M$ of the polygon $P$, and the generators of
its  coefficient semifield are in bijection with the boundary edges
$T_{n+1},\ldots,T_{2n+3}$  of $P$.
To be more precise, the coefficient semifield is the tropical
semifield $\textup{Trop}(x_{n+1},x_{n+2},\ldots,x_{2n+3})$, which is
a free abelian group, written multiplicatively, with generators
$x_{n+1},\ldots,x_{2n+3}$, and with auxiliary addition $\oplus$ given
by 
\[ \prod_j x_j^{a_j} \oplus \prod_j x_j^{b_j}= \prod_j x_j^{\textup{min}(a_j,b_j)}.
\] 
Clusters are in bijection with triangulations of $P$. Given a
triangulation $T=\{T_1,\ldots,T_n,T_{n+1},\ldots,T_{2n+3}\}$ let
$\mathbf{x}_T=\{x_1,\ldots,x_n\}$ be the corresponding cluster, where
we use the notation $x_i=x_{T_i} $ for short. The mutation in
direction $k$ is
described as follows:  $T_k$ is  a diagonal in a unique  quadrilateral
in $T$. Let $T_{k'}$ be the other diagonal in that quadrilateral. Then
the mutation in direction $k$ of $T$ is the triangulation obtained
from $T$ by replacing ${T_k}$ by ${T_{k'}}$. The corresponding exchange
relation is $x_kx_{k'}=x_ax_c+x_bx_d$, where $T_a,T_c$ are two
opposite sides  of the quadrilateral and  $T_b,T_d$ are the  other two
opposite sides.

Let us choose the initial seed $(\mathbf{x},\mathbf{p},B)$ 
consisting of a cluster $\mathbf{x}$, a coefficient vector $
\mathbf{p}=(p_1^+,p_1^-,p_2^+,p_2^-,\ldots,p_n^+,p_n^-)\in
\left(\textup{Trop}(x_{n+1},x_{n+2},\ldots,x_{2n+3})\right) ^{2n}$ 
and a $(2n+3)\times n$  matrix $B=(b_{ij})$ as follows:
Taking the initial cluster  to be the snake
triangulation of \cite[section 12]{FZ2}, the initial matrix
$B=(b_{ij})$ is given by the conditions 
$b_{ii}=0, b_{ij}\in \{-1,0,1\} $ and $b_{ij}=1$ (respectively
$b_{ij}=-1$) if and only
if $i\ne j$, the edges $T_i$ and $T_j$  bound the same triangle
and the sense of rotation from $T_i$  to $T_j$ is counterclockwise
(respectively clockwise).
The corresponding initial coefficient vector is given by
\[p_j^+=\prod_{i\,\ge\, n+1 \,: \, b_{ij}=1 } x_j \qquad \textup{and} \qquad  
p_j^-=\prod_{i\,\ge\, n+1 \,: \, b_{ij}=-1} x_j.\]
An example of rank $3$ is given in Figure \ref{fig 1}.
\begin{figure}
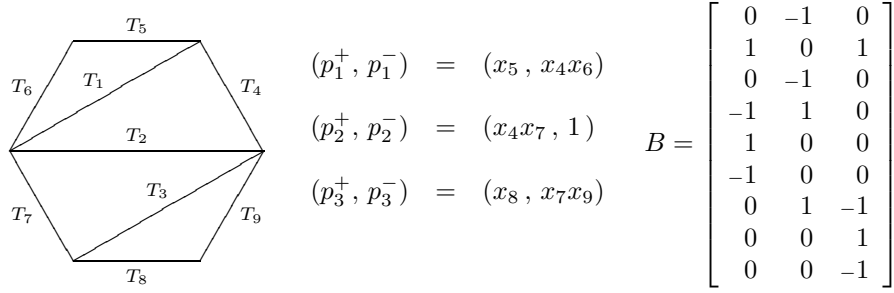

\def\alphanum{\ifcase \xypolynode \or 4 \or 5\or 6 \or 7\or 8\or 9\fi}
$$\begin{array}{ccc}
\xy/r4pc/: {\xypolygon6"A"{~<<{@{}}~><{@{-}}
~>>{_{T_{\alphanum}}}}},
\POS"A2" \ar@{-}_{T_1} "A4",
\POS"A4" \ar@{-}^{T_2} "A1",
\POS"A5" \ar@{-}^{T_3} "A1",
\endxy &
\begin{array}{rcl}
(p_1^+,\,p_1^-)&=&(x_5\,,\,x_4x_6)\\ \\
(p_2^+,\,p_2^-)&=&(x_4x_7\,,\,1\,)\\ \\
(p_3^+,\,p_3^-)&=&(x_8\,,\,x_7x_9)\\ \\
\end{array}  
&B=\left[\! 
\begin{array}{rrr} 
0&{\scriptstyle -} 1&0\\ 1&0&1\\ 0&{\scriptstyle -}1&0\\ {\scriptstyle -}1&1&0\\1&0&0\\ {\scriptstyle -}1&0&0\\ 0&1&{\scriptstyle -}1\\ 0&0&1\\0&0&{\scriptstyle -}1\\
\end{array}  \right]
\end{array}  $$
\caption{Snake triangulation, initial coefficients and initial
  matrix in rank $3$}\label{fig 1}
\end{figure}   
\end{subsection} 

\begin{subsection}{$T$-paths}\label{sect 1.2}

Let $T=\{T_1,T_2,\ldots,T_n,T_{n+1},\ldots,T_{2n+3}\}$ be a
 triangulation of the polygon $P$, where $T_1,\ldots,T_n$ are
 diagonals and $T_{n+1},\ldots,T_{2n+3}$ are boundary edges. 
 Let $a$ and $b$ be two
non-adjacent vertices on the boundary and let $M_{a,b}$ be the
diagonal that connects $a$ and $b$.

\begin{definition}\label{Tpath}
A \emph{$T$-path} $\za$ from $a$ to $b$ is  a sequence 
\[ \za = (a_0,a_1,\ldots,a_{\ell(\za)}\mid
i_1,i_2,\ldots,i_{\ell(\za)})\]
such that 
\begin{itemize}
\item[\textup{(T1)}] $a\!=\!a_0,a_1,\ldots,a_{\ell(\za)}\!=\!b$ are vertices of $P$.
\item[\textup{(T2)}] $i_k\in \{1,2,\ldots,{2n+3}\}$ such
  that $T_{i_k}$ connects the vertices $a_{i_{k-1}}$ and $a_{i_k}$
  for each $k=1,2,\ldots,\ell(\za)$. 
\item[\textup{(T3)}] $i_j\ne i_k$ if $j\ne k$.
\item[\textup{(T4)}] $\ell(\za)$ is odd.
\item[\textup{(T5)}] $T_{i_{k}}$ crosses $M_{a,b}$ if    $k$ is even.
\item[\textup{(T6)}] If  $j<k$ and both $T_{i_j}$ and $T_{i_{k}}$ cross
  $M_{a,b}$ then the crossing point of $T_{i_j}$  
and $M_{a,b}$ is closer to the vertex $ a $   
than the crossing point of $T_{i_{k}}$ and $M_{a,b}$.

\end{itemize}    

\end{definition}

Thus, a $T$-path from $a$ to $b$ is a path on the edges of the
triangulation $T$, that does not use any edge twice,
 whose crossing points with $M_{a,b}$ are progressing from $a$
 towards $b$, and, when 
classifying the edges into even 
and odd edges according to their order of appearance, then
every even edge is crossing $M_{a,b}$.

To any $T$-path $\za = (a_0,a_1,\ldots,a_{\ell(\za)}\mid
i_1,i_2,\ldots,i_{\ell(\za)})$, 
we associate an element $x(\za)$ in the cluster algebra by

\begin{equation}\label{**}
x(\za)=
{\prod_{k \textup{ odd}} x_{i_{k}}}
\
{\prod_{k \textup{ even}}x_{i_{k}}^{-1}} . 
\end{equation}

\begin{definition}
Let $\mathcal{P}_T(a,b)$ denote the
set of $T$-paths  from $a$ to $b$.
\end{definition}     

\begin{lem}\label{Lemma 1}
 Let $\alpha,\alpha'\in \mathcal{P}_T(a,b)$ with $\alpha \ne \alpha'$. Then 
$x(\alpha)\ne x(\alpha')$.
\end{lem}
\begin{pf}
Suppose $x(\alpha)= x(\alpha')$. Then (T3) implies that the set of 
even edges and the set of odd edges are the same in $\alpha$ and $\alpha'$. 
From conditions (T5) and (T6) it follows that the order of the even
edges is the same,  
hence the order of all edges is the same, whence $\alpha=\alpha'$.
\qed
\end{pf}

\end{subsection}

\begin{subsection}{Expansion formula}\label{sect main}
The following theorem is our main result. 

\begin{thm}\label{thm 1} Let $a$ and $b$ be two non-adjacent vertices of
  $P$, let $M=M_{a,b}$ be the diagonal connecting $a$ and $b$ and let
  $x_M$ be the corresponding cluster variable. Then

\begin{equation}\label{*}
  x_{M}=\sum_{\za\in\mathcal{P}_T(a,b)} x(\za).
\end{equation}
\end{thm}  

\begin{rem}\label{rem 1} \textup{Because of conditions \textup{(T3)} and \textup{(T5)}, 
each $x(\za)$ is a reduced fraction whose denominator is a product of
cluster variables.}
\end{rem}
\begin{rem}\label{rem 2}\textup{By Lemma \ref{Lemma 1},
    each term   in the sum of equation 
  \textup{(\ref{*})}  appears with multiplicity one.}
\end{rem} 
\begin{rem}\label{rem 3} \textup{ Formula \textup{(\ref{*})} also applies to
type $A$  cluster algebras with trivial coefficients by setting
$x_t=1$ for $t=n+1,n+2, \ldots,2n+3$.}
\end{rem} 
The proof of Theorem \ref{thm 1} will be given in section \ref{sect proof}. To
illustrate the statement, we give an example here.

\begin{example}
The following figure shows a triangulation $T=\{T_1,\ldots,T_{13}\}$ and
a (dotted) diagonal $M=M_{a,b}$. Next to it is a complete list
 of elements of $\mathcal{P}_T(a,b)$. 
\[\begin{array}{cc}
\def\alphanum{\ifcase \xypolynode \or 6 \or 7\or 8 \or 9\or 10\or
  11\or 12\or 13\fi}
\xy/r6pc/: {\xypolygon8"A"{~<<{@{}}~><{@{-}}
~>>{_{T_{\alphanum}}}}},
\POS"A2" \ar@{-}^(0.7){T_1} "A4",
\POS"A4" \ar@{-}_(0.5){T_2} "A6",
\POS"A6" \ar@{-}^(0.3){T_3} "A2",
\POS"A2" \ar@{-}^{T_4} "A8",
\POS"A6" \ar@{-}^(0.3){T_5} "A8",
\POS"A3" \ar@{.}^(0.6){M} "A7",
\POS"A3"\drop{\begin{array}{c}a\ \\ \\ \end{array}  }
\POS"A7"\drop{\begin{array}{c} \\ \ b \end{array}  } 
\POS"A2"\drop{\begin{array}{c}\  f\\ \\ \end{array}  } 
\POS"A4"\drop{\begin{array}{c} c\quad \\ \end{array}  }
  \POS"A6"\drop{\begin{array}{c} \\d \ \ \end{array}  }
  \POS"A8"\drop{\begin{array}{c} \quad e \end{array}  }  
\endxy
&

\begin{array}{c}
 (a,f,d,b\mid 7,3,11)\, 
\\(a,f,c,d,e,b\mid 7,1,2,5,12)\\
(a,c,f,e,d,b\mid 8,1,4,5,11)
\\(a,c,f,d,e,b\mid 8,1,3,5,12)\\
(a,f,c,d,f,e,d,b \mid 7,1,2,3,4,5,11) 
\end{array} 
\end{array}  \]
Theorem \ref{thm
  1} thus implies that 
\[ x_M=
\frac{x_7x_{11}}{x_3} +
\frac{x_7x_2 x_{12}}{x_1x_5} +
\frac{x_8x_4 x_{11}}{x_1x_5} +
\frac{x_8x_3 x_{12}}{x_1x_5} +
\frac{x_7x_2x_4 x_{11}}{x_1x_3x_5} .
\]
\end{example}

\end{subsection} 

\begin{subsection}{Positivity}
In the case of the Ptolemy algebra, the following positivity
conjecture of \cite{FZ1} is a direct 
consequence of Theorem \ref{thm 1} and Remarks \ref{rem 1} and
\ref{rem 2}.

\begin{cor}\label{cor 1}
Let $x$ be any cluster variable in the Ptolemy cluster algebra and let
$\{x_1,\ldots,x_n\}$ be any cluster. Let
\[x=\frac{f(x_1,\ldots,x_n,x_{n+1},\ldots,x_{2n+3})}{x_1^{d_1}\ldots x_n^{d_n}}\] be the
expansion of $x$ in the cluster $\{x_1,\ldots,x_n\}$, where $f$ is a
polynomial which is not divisible by any of the $x_1,\ldots,x_n$.
Then the coefficients of $f$ are  either $0$ or $1$, thus, in
particular, they are non-negative integers.
\end{cor}

\end{subsection}
\end{section}

\begin{section}{Proof of Theorem \ref{thm 1}}\label{sect proof}
This section is devoted to the proof of our main result. Let $T$ be a
triangulation of the polygon $P$, let $a$ and $b$ be two vertices of
$P$ and $M=M_{a,b}$ the diagonal connecting $a$ and $b$. Suppose 
that $M\notin T$.
Among all diagonals of $T$ that cross
$M$, there is a unique one, $T_{i_0}$, such that its intersection point
with $M$ is the closest possible to the vertex $a$. 
Then there is a
unique triangle in $T$ having $T_{i_0}$ as one side and the vertex $a$
as third point. Denote the other two sides of this triangle by
$T_{i_1} $ and $T_{i_1'}$ and let $c$ be the common endpoint of
$T_{i'_1}$ and $T_{i_0}$, and  $d$  the common endpoint of
$T_{i_1}$ and $T_{i_0}$ (see Figure \ref{fig 3}). Note that $T_{i_1},
  T_{i_1'}$ may be boundary edges.
\begin{figure}
\begin{center}
\includegraphics{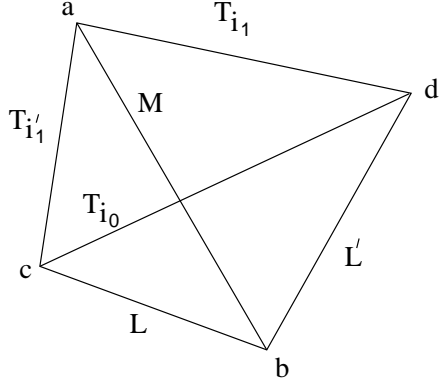}
\caption{Proof of Theorem \ref{thm 1}}\label{fig 3}
\end{center}
\end{figure}   
Now consider the unique quadrilateral  in which $M$ and
$T_{i_0}$ are the diagonals. Two of its sides are $T_{i_1}$ and
$T_{i_1'}$. Denote the other two sides by $L$ and
$L'$ in such a way that $L$ is the side opposite to $T_{i_1}$
(see Figure \ref{fig 3}). We will keep this setup for the rest of this
section.

\begin{lem}\label{lem 0}
\begin{itemize}
\item[(a)] If $T_i\in T$ crosses $L$ (respectively $L'$), then $T_i$
  crosses $M$.
\item[(b)]  If $T_i\in T$ is adjacent to $a$, then $T_i$ does not
  cross $L$ nor $L'$.
\item[(c)]  If $T_i\in T$ crosses $M$ and does not cross $L$
  (respectively $L'$), then $T_i$ is adjacent to $c$ (respectively $d$).
\end{itemize}   
\end{lem}  

\begin{pf}
This follows directly from the construction.
\qed
\end{pf}  

Let $\mathcal{P}_T(a,b)_{j}$ denote the subset of $\mathcal{P}_T(a,b)
$  of all 
$T$-paths $\za$ that start with the edge $T_j$ and let
$\mathcal{P}_T(a,b)_{-j}$ be the subset of $\mathcal{P}_T(a,b) $
 of all 
$T$-paths $\za$ that do not contain  the edge $T_j$.
Similarly, let $\mathcal{P}_T(a,b)_{jk}$ denote the subset of
$\mathcal{P}_T(a,b)_j $   of all
$T$-paths $\za$ that start with the edges $T_jT_k$ and let
$\mathcal{P}_T(a,b)_{j,-k}$ be the subset of $\mathcal{P}_T(a,b)_j $ of all  
$T$-paths $\za$ that start with the edge $T_j$ and do not contain
the edge $T_k$.

\begin{lem}\label{lem 5} We have
$\mathcal{P}_T(a,b)=\mathcal{P}_T(a,b)_{i_1}\sqcup
\mathcal{P}_T(a,b)_{i_1'}$. 
\end{lem}

\begin{pf}
Let $\za=(a_0,a_1,\ldots,a_{\ell(\za)}\mid
{j_1},{j_2},\ldots,{j_{\ell(\za)}})$ be an
arbitrary element of 
$\mathcal{P}_T(a,b)$. 
We know that $T_{j_1}\in T$ and that one of its
  endpoints is the vertex $a$. Moreover, $T_{j_2}$ crosses $M$, by
  property (T5) of $T$-paths.
Let $P_a$ and $P_b$ denote the
  two pieces of the polygon  obtained by cutting along
  $T_{i_0}$, where the vertex $a$ lies in the piece $P_a$, while the
  vertex $b$ lies in $P_b$.
  The piece $P_a$ also contains the
  diagonal $T_{j_1}$ whereas the piece $P_b$ contains all the diagonals of
  $T$ that cross $M$ and, therefore, $P_b$ contains $T_{j_2}$.
  Consequently, the diagonals $T_{j_1},T_{j_2}$ and $T_{i_0}$ share
  one endpoint, which is either $c$ or $d$. 
There are therefore  two possibilities: either $T_{j_1}=T_{i_1}$ or
$T_{j_1}=T_{i_1'} $.
\qed
\end{pf}

\begin{lem}\label{cor 2}
\begin{itemize}
\item[(a)] $\mathcal{P}_T(a,b)_{i_1}=\mathcal{P}_T(a,b)_{i_1i_0}
  \sqcup \mathcal{P}_T(a,b)_{i_1,-i_0}$
\item[(b)] $\mathcal{P}_T(a,b)_{i'_1}=\mathcal{P}_T(a,b)_{i'_1i_0}
  \sqcup \mathcal{P}_T(a,b)_{i'_1,-i_0}$
\end{itemize}   
\end{lem}   %

\begin{pf} By constuction, $T_{i_0}$ is the edge of the triangulation $T$ such
  that its crossing point with $M$ is the closest possible to the
  vertex $a$. Hence, the result follows from condition (T6).
\qed
\end{pf}

Now let $\zg=(c_0,\ldots,c_{\ell(\zg)}\mid j_1,\ldots,j_{\ell(\zg)})$ be
any $T$-path from $c$ to $b$.
Suppose first that $j_1=i_0$, thus $c_1=d$. 
In this case, let $f(\zg)$ be the path
obtained from 
$\zg$ by replacing the first edge $i_0 $ by $i_1$, that is 
\[f(\zg)=(a,c_1,\ldots,c_{\ell(\zg)}\mid i_1,j_2,\ldots,j_{\ell(\zg)}).\]
Suppose now that $j_k\ne i_0$ for all $k$. In this case, let $g(\zg)$ be the
composition of the paths $(a,d,c\mid i_1,i_0)$ and
$\zg$, that is 
\[g(\zg)=(a,d,c_0,c_1,\ldots,c_{\ell(\zg)}\mid
i_1,i_0,j_1,j_2,\ldots,j_{\ell(\zg)}).\] 
Let us check that $f(\zg)$ and $g(\zg)$ are elements of
$\mathcal{P}_T(a,b)$.
 Indeed, the properties (T1),(T2) and (T4) are immediate and (T5)
follows from Lemma \ref{lem 0}(a). In order to show (T3), we need to
prove that $i_1\ne j_k$ for all $k$; but this follows from Lemma
\ref{lem 0}(b) and the fact that $T_{i_1}$ is adjacent to $a$.
 Finally, (T6) holds since the crossing point of  $T_{i_0}$ 
and $M_{a,b}$ is the closest possible to $a$.
We have the following lemma.
\begin{lem}\label{lem f}
The maps $f$ and $g$  induce  bijections 
\[
\begin{array}{ccccccccc}
 f:\mathcal{P}_T(c,b)_{i_0}\to
 \mathcal{P}_T(a,b)_{i_1,-i_0}&\quad and \quad&
 g:\mathcal{P}_T(c,b)_{-i_0}\to \mathcal{P}_T(a,b)_{i_1i_0},
\end{array}  \]
and 
\begin{equation}\label{eq 2}
  x(f(\zg))=\frac{x_{i_1}}{x_{i_0}} \, x(\zg) \quad \textup{\it and}\quad
 x(g(\zg))=\frac{x_{i_1}}{x_{i_0}} \, x(\zg) 
\end{equation}
\end{lem}  
\begin{pf} 
The formulas (\ref{eq 2}) follow directly from the definitions of $f$
and $g$. These formulas together with Lemma \ref{Lemma 1} imply the
injectivity of $f$ and $g$. 
To show the surjectivity of $f$, suppose that 
$\mathcal{P}_T(a,b)_{i_1,-i_0}$ is not empty and
let $\za\in \mathcal{P}_T(a,b)_{i_1,-i_0}$ be an arbitrary
element. Say 
\[\za=(a,d,a_2,a_3,\ldots,a_{\ell(\alpha)} \mid
i_1,j_2,j_3,j_4,\ldots,j_{\ell(\alpha)}).\]
We need to show that the path
\[\zg=(c,d,a_2,a_3,\ldots,a_{\ell(\alpha)} \mid 
i_0,j_2,j_3,j_4,\ldots,j_{\ell(\alpha)})\]
is an element of $\mathcal{P}_T(c,b)_{i_0}$.
Conditions (T1),(T2),(T4)  hold since $\za\in
\mathcal{P}_T(a,b)$, condition (T3) holds because the path $\za$
does not contain the edge $T_{i_0}$ and condition (T6) holds because
of Lemma \ref{lem 0}(a). 
Moreover, since 
$\mathcal{P}_T(a,b)_{i_1,-i_0}$ is not empty,  there exists a
  diagonal in $T\setminus\{T_{i_0}\}$ which is adjacent to $d$ and
  crosses $M$. Since $T$ is a triangulation, it follows that any
  diagonal in  $T\setminus\{T_{i_0}\}$ that crosses $M$ also crosses
  $L$. Thus $\zg$ satisfies condition (T5), because  $\za\in
\mathcal{P}_T(a,b)$. 
This shows that $\zg\in \mathcal{P}_T(c,b)$, and since $\zg$ starts
with the edge $i_0$, we have  $\zg\in \mathcal{P}_T(c,b)_{i_0}$.
 Hence  $f$ is surjective.

It remains to show that $g$ is surjective. Let $\alpha \in
\mathcal{P}_T(a,b)_{i_1i_0}$ be arbitrary. Say 
\[\alpha=(a,d,c,a_3,\ldots,a_{\ell(\alpha)} \mid 
i_1,i_0,j_3,j_4,\ldots,j_{\ell(\alpha)}).\]
We have to show that
\[\gamma=(c,a_3,\ldots,a_{\ell(\alpha)} \mid j_3,j_4,\ldots,j_{\ell(\alpha)})
\in \mathcal{P}_T(c,b).\]
Conditions (T1)-(T4) hold for $\gamma$ because $\alpha \in
\mathcal{P}_T(a,b)$ and condition (T6) holds because
of Lemma \ref{lem 0}(a).  Let us show (T5).  
We need to show that any even edge of $\gamma$ crosses $L$. Since
$\alpha\in\mathcal{P}_T(a,b)$, we know that every even edge of $\gamma$ crosses $M$. 
Thus by Lemma \ref{lem 0}(c), if there is an even edge of $\gamma$ that does not cross $L$ then 
this edge has to be adjacent to $c$.
Since $\gamma$ starts at $c$, its first even edge $T_{j_4}$ cannot  be 
adjacent to $c$, and thus $T_{j_4}$ crosses both $M$ and $L$. Then, since 
$\alpha $ satisfies (T6), every even edge of $\gamma$ crosses $M$ and $L$. 
This shows (T5).
Hence $\gamma\in \mathcal{P}_T(c,b)$ and $g$ is surjective.
\qed
\end{pf}

\begin{lem}\label{lem sum} We have
\[
\begin{array}{lcrclc}
(a)&\hspace*{45pt}& \displaystyle\sum_{\zg\in \mathcal{P}_T(c,b)  }
  x(\zg)\ \frac{x_{i_1}}{x_{i_0}} &=&   
\displaystyle \sum_{\za\in
 \mathcal{P}_T(a,b)_{i_1}  } x(\za)&\hspace*{50pt} \\ \\ 
(b)&& \displaystyle\sum_{\zg\in \mathcal{P}_T(d,b)  } x(\zg) \
  \frac{x_{i'_1}}{x_{i_0}} &=& 
\displaystyle\sum_{\za\in 
 \mathcal{P}_T(a,b)_{i'_1}  } x(\za).
\end{array}  
\]
\end{lem}  

\begin{pf}
The first statement follows from  Lemma \ref{cor 2}(a), Lemma
\ref{lem f} and the fact that 
 $\mathcal{P}_T(c,b) = \mathcal{P}_T(c,b)_{i_0}\sqcup
\mathcal{P}_T(c,b)_{-i_0}$. 
The second statement follows by symmetry.
\qed
\end{pf}

Proof of Theorem \ref{thm 1}.
For any $T_i\in T$, let $e(T_i,M) \in\{0,1\}$ be the number of
crossings of the diagonals $T_i$ and $M$. Then the total number
of crossings between $M$ and $T$ is $e(T,M)=\sum_{T_i\in
  T} e(T_i,M)$. 

We prove the theorem by induction on $e(T,M)$.
If $e(T,M) =0$, then $M=T_i\in T$ for some
$i\in\{1,\ldots,n\}$. In this case, no element of $T$ crosses
$M$ and, by condition (T5), the set
$\mathcal{P}_T(a,b)$ contains exactly one element: $(a,b\mid i)$. Thus 
\[\sum_{\za\in \mathcal{P}_T(a,b) }x(\za) = x(a,b\mid i) = x_i =
x_{M}.\]

Suppose now that $e(T,M)\ge 1$. 
As before, consider the unique quadrilateral in $T$ in which $M$ and
$T_{i_0}$ are the diagonals (see Figure \ref{fig 3}).
Thus, in the cluster algebra, we have the following  exchange relation 
\begin{equation}\label{exch}
  x_M\,x_{i_0} = x_{i_1}\,x_L +x_{i_1'}\,x_{L'}. \end{equation} 
Moreover, any diagonal in $T$ that crosses $L$ (respectively $L'$) also
crosses $M$, by Lemma \ref{lem 0}(a), and, moreover, $T_{i_0}$ crosses
$M$ but crosses neither $L$ nor $L'$. Thus $e(T,L)<e(T,M)$ and $
e(T,L')<e(T,M)$, and by 
induction hypothesis
\[x_L=\sum_{\zg\in \mathcal{P}_T(c,b)} x(\zg)
\qquad \textup{and} \qquad
x_{L'}=\sum_{\zg\in \mathcal{P}_T(d,b)} x(\zg).
\]
Therefore, we can write the exchange relation (\ref{exch}) as 
\begin{eqnarray}\nonumber
 x_M &=&\sum_{\zg\in  \mathcal{P}_T(c,b)} x(\zg) \ \frac{x_{i_1}}{x_{i_0}}
\quad+\sum_{\zg\in \mathcal{P}_T(d,b)} x(\zg)\ \frac{x_{i'_1}}{x_{i_0}}.
\end{eqnarray}  
The theorem now follows from Lemma \ref{lem sum} and Lemma \ref{lem 5}.
\qed

\end{section} 

{} 

\bigskip 
\bigskip 

\noindent 
Department of  Mathematics and    Statistics\\
University of   Massachusetts at Amherst\\
Amherst, MA 01003--9305, USA\\ 
E-mail address: {\tt schiffler@math.umass.edu}

\end{document}